\documentclass{article}
\usepackage{authblk}
\usepackage{graphicx}
\usepackage{epstopdf}

\begin{document}

\title{Lie symmetries of (1+2) nonautonomous evolution equations in
Financial Mathematics}
\author[1]{A Paliathanasis\thanks{%
paliathanasis@na.infn.it}}
\author[2]{RM Morris\thanks{%
rmcalc85@gmail.com}}
\author[2,3,4]{ PGL Leach\thanks{%
leach@ucy.ac.cy}}

\affil[1]{Instituto de Ciencias F\'{\i}sicas y Matem\'{a}ticas, Universidad Austral de
Chile, Valdivia, Chile}
\affil[2]{Department of Mathematics and
Institute of Systems Science, Research and Postgraduate Support, Durban
University of Technology, PO Box 1334, Durban 4000, Republic of South Africa}
\affil[3]{School of Mathematics, Statistics and Computer Science, University
of KwaZulu-Natal, Private Bag X54001, Durban 4000, Republic of South Africa}
\affil[4]{Department of Mathematics and Statistics, University of Cyprus,
Lefkosia 1678, Cyprus}

\renewcommand\Authands{ and }
\maketitle

\begin{abstract}
We analyse two classes of $(1+2)$ evolution equations which are of special
interest in Financial Mathematics, namely the Two-dimensional Black-Scholes
Equation and the equation for the Two-factor Commodities Problem. Our
approach is that of Lie Symmetry Analysis. We study these equations for the case in which they are autonomous
 and for the case in which the parameters of the equations are
unspecified functions of time. For the
autonomous Black-Scholes Equation we find that the symmetry is maximal and
so the equation is reducible to the $(1+2)$ Classical Heat Equation. This is
not the case for the nonautonomous equation for which the number of
symmetries is submaximal. In the case of the two-factor equation the number
of symmetries is submaximal in both autonomous and nonautonomous cases. When
the solution symmetries are used to reduce each equation to a $(1+1)$
equation, the resulting equation is of maximal symmetry and so equivalent to
the $(1+1)$ Classical Heat Equation.
\end{abstract}


\section{Introduction}

In the early seventies F Black and M Scholes \cite{bsch01,bsch02} and,
independently, R Merton \cite{Merton} introduced a mathematical model for
the pricing of European options. \ The Black-Scholes-Merton (BS) Model is described by an
$(1+1)$ evolution equation. The mathematical expression of the BS equation
is
\begin{equation}
\frac{1}{2}\sigma ^{2}S^{2}u_{,SS}+rSu_{,S}-ru+u_{,t}=0,  \label{bs.01}
\end{equation}%
in which $t$ is time, $S$ is the current value of the underlying asset, for
example a stock price, $r$ is the rate of return on a safe investment, such
as government bonds and $u=u\left( t,S\right) $ is the value of the option.  The solution of (\ref{bs.01} is
subject to the satisfaction of the terminal condition $u\left( T,S\right) =U$%
, when $t=T$.

For the prices of commodities, E Schwartz \cite{schwartz} proposed three
models which study the stochastic behaviour of the prices of commodities
that take into account several aspects of possible influence on the prices.
\ In the simplest model he assumed that the logarithm of the spot price
followed a mean-reversion process of Ornstein-Uhlenbeck type. This is termed
the one-factor model. The one-factor model is described by the equation
\begin{equation}
\frac{1}{2}\sigma ^{2}S^{2}F_{,SS}+\kappa \left( \mu -\lambda -\log S\right)
Su_{,S}-F_{,t}=0,  \label{bs.02}
\end{equation}%
where $\kappa >0$ measures the degree of reversion to the long-run mean
log price, $\lambda $ is the market price of risk, $\mu $ is the drift rate
of $S$ and $F=F\left( t,S\right) $ is the current value of the futures
contract. The solution of (\ref{bs.02}) satisfies the initial condition $%
F\left( 0,S\right) =S.$

The BS equation (\ref{bs.01}) and the one-factor equation (\ref{bs.02}) are
of the same equivalence class as the Schr\"odinger equation and the Heat
diffusion equation. All four equations model random phenomena of different
contexts. The two first are in financial mathematics, the third in quantum
physics and the fourth in dispersion.

It has been proven that all four equations are maximally symmetric and
invariant under the same group of invariant transformations of dimension $%
5+1+\infty $ which span the Lie algebra $\left\{ sl\left( 2,R\right) \oplus
_{s}W_{3}\right\} \oplus _{s}\infty A_{1}$, where $W_3$ is a representation of the three-dimensional Weyl--Heisenberg Group,
(in the Mubarakzyanov Classification Scheme \cite%
{Morozov58a,Mubarakzyanov63a,Mubarakzyanov63b,Mubarakzyanov63c} this is $\left\{ A_{3,8}\right) \oplus
_{s} A_{3,1}\} \oplus _{s}\infty A_{1}$).  This means
that there exists a point transformation which transforms one equation to
another. The Lie symmetries of the BS equation (\ref{bs.01}) have been found
in \cite{ibrag}, whereas the Lie symmetries of the one-factor model (\ref%
{bs.02}) were found in \cite{Leach08}. \

The parameters of the models (\ref{bs.01}) and (\ref{bs.02}) ) are generally
assumed to be constant. However, in real problems they may vary with time if
the time-span of the model is sufficiently long. In \cite{Tamizhmani} it has
been shown that, when the parameters $\sigma $, and $r$ of the BS equation
are time-dependent, \textit{ie}, $\sigma =\sigma \left( t\right) $ and $%
r=r\left( t\right) $, the time-dependent BS equation is invariant under the
same group of invariant transformations as that of the ``static'' BS equation.
The same result has been found for the time-dependent one-factor model of
commodities \cite{onefactor}. Hence the autonomous and the nonautonomous
equations (\ref{bs.01}) and (\ref{bs.02}) are maximally symmetric and
equivalent under point transformations. \

In Classical Mechanics the slowly lengthening pendulum with equation of
motion in the linear approximation,%
\begin{equation}
\ddot{x}+\omega ^{2}\left( t\right) x=0,  \label{bs.03}
\end{equation}%
in which the time dependence in the `spring constant' is due to the length
of the pendulum's string increasing slowly \cite{werner}, admits the
conservation law \cite{Lewis1,Lewis2} (note that the case of a slowly
shortening pendulum is quite different \cite{slowly})%
\begin{equation}
I=\frac{1}{2}\left\{ \left( \rho \dot{x}-\dot{\rho}x\right) +\left( \frac{x}{%
\rho }\right) ^{2}\right\} ,  \label{bs.04}
\end{equation}%
where $\rho =\rho \left( t\right) $, is a solution of the second-order
differential equation%
\begin{equation}
\ddot{\rho}+\omega ^{2}\left( t\right) \rho =\frac{1}{\rho ^{3}}.
\label{bs.05}
\end{equation}
This result is independent of the rate of change of the length of the
pendulum.

The latter equation is the well-known Ermakov-Pinney equation \cite{Ermakov}%
. The solution was given by Pinney in \cite{Pinney} and it is
\begin{equation}
\rho \left( t\right) =\sqrt{A\upsilon _{1}^{2}+2B\upsilon _{1}\upsilon
_{2}+C\upsilon _{2}^{2}}  \label{bs.06}
\end{equation}%
subject to a constraint on the three constants, $A$, $B$ and $C$. Functions $%
\upsilon _{1}\left( t\right) ,~\upsilon _{2}\left( t\right) ,$ are two
linearly independent solutions of (\ref{bs.03}) .

Equation (\ref{bs.03}) is invariant under the action of
the group invariant transformations in which the generators of the infinitesimal transformations
form the $sl\left( 3,R\right) $ algebra. This is the Lie algebra admitted by
the harmonic oscillator, $\omega \left( t\right) =\omega _{0}$, and the
equation of the free particle,~$\omega \left( t\right) =0$ \cite%
{Lut,Leach1980,LeachA08}. The transformation which connects the nonautonomous
linear equation (\ref{bs.03}) with the autonomous oscillator is a
time-dependent linear canonical transformation of the form
\begin{equation}
Q=\frac{x}{\dot{x}}~\ ,~P=\rho \dot{x}-\dot{\rho}x~,~T=\int^{t}\rho
^{-2}\left( \eta \right) d\eta,  \label{bs.07}
\end{equation}%
where $\rho $ is given by (\ref{bs.06}).

The connection of the number of symmetries of the corresponding
Schr\"'odinger Equation with the Noether point symmetries of the classical
Lagrangian \cite{Leach1,And1} was seen to extend to the time-dependent case
\cite{Maharaj2010} and, indeed, be seen to be the same as the equivalent
autonomous systems \cite{Leach2005} and in the case of maximal symmetry is $%
\left\{ sl\left( 2,R\right) \oplus _{s}W_{3}\right\} \oplus _{s}\infty A_{1}$
which is that of the $(1 + 1)$ classical heat equation.

In this context we wish to see what happens when we pass from an autonomous $%
(1+2)$ evolution equation to the corresponding nonautonomous case. For that
we study the Lie symmetries of the nonautonomous models\ of: (a) the
two-factor model of commodities and (b) the two-dimensional BS equation. \

We find that, for the two-factor model, the autonomous and the nonautonomous
equations are invariant under the same group of invariant transformations $%
\left\{ A_{1}\oplus _{s}W_{5}\right\} \oplus _{s}\infty A_{1}$. However,
that it is not true for the two-dimensional BS equation. The reason for that
is that the Lie symmetries of the two-factor model follow from the
translation group of the two-dimensional Euclidian space (except the
homogeneous and the infinite number of solution symmetries). The translation group generates
Lie symmetries for both the autonomous system and for the nonautonomous
system.

On the other hand the autonomous two-dimensional BS equation is maximally
symmetric, \textit{ie}, it admits nine Lie symmetries plus the infinite number of solution
symmetries, which form the $\left\{ \left\{ sl\left( 2,R\right)
\oplus _{s}so\left( 2\right) \right\} \oplus _{s}W_{5}\right\} \oplus
_{s}\infty A_{1}~$Lie algebra.$~$This result completes the analysis of \cite%
{consta} in which they found that the two-dimensional\ BS equation admits
seven Lie point symmetries plus the $\infty A_{1}$.

The nonautonomous two-dimensional BS equation is invariant under the Lie
algebra \newline
$\left\{ \left\{ A_{1}\oplus _{s}so\left( 2\right) \right\} \oplus
_{s}W_{5}\right\} \oplus _{s}\infty A_{1}$, that is, the $sl\left(
2,R\right) $ subalgebra is lost. The reason for that is that the Lie symmetries
of the autonomous two-dimensional BS equation arise from the homothetic
algebra of the two-dimensional Euclidian space which defines the Laplace
operator of the evolution equation and, when the parameters in the second
derivatives are not constants, the homothetic algebra of the Euclidian
space does not generate Lie symmetries. Moreover, in the case for which the
parameters of the second derivatives are time-indepedent, the
two-dimensional BS equation is maximally symmetric, \textit{ie}, it is
invariant under the same group of point transformations as the (1+2)
autonomous BS and Heat conduction equations.

The plan of the paper is as follows.
In Section \ref{twofactor} we study the Lie symmetries of the two-factor
model of commodities for the autonomous and nonautonomous cases. We show
that in both cases the two-factor model is invariant under the $\left\{
A_{1}\oplus _{s}W_{5}\right\} \oplus _{s}\infty A_{1}$ Lie algebra. The Lie
symmetries of the two-dimensional BS equation, the autonomous and the
nonautonomous, are studied in Section \ref{2DBS}. Finally in Section \ref{con}
we give some applications and we draw our conclusions.

\section{The two-factor model of commodities}

\label{twofactor}

The two-factor model adds to the spot price, $S$, of (\ref{bs.02}) the
instantaneous convenience yield, $\delta $, which may be interpreted as the
flow of services accruing to the holder of the spot commodity but not to the
owner of a futures contract. The evolution partial differential equation for
this model is
\begin{equation}
\frac{1}{2}\sigma _{1}^{2}S^{2}F_{,SS}+\rho \sigma _{1}\sigma
_{2}F_{,S\delta }+\frac{1}{2}\sigma _{2}^{2}F_{,\delta \delta }+\left(
r-\delta \right) SF_{,S}+\left( \kappa \left( \alpha -\delta \right)
-\lambda \right) F_{,\delta }-F_{,t}=0  \label{bs.08}
\end{equation}%
for which the terminal condition is now $F\left( 0,S,\delta \right) =S.$

Equation (\ref{bs.08}) is an $(1+2)$ evolution equation and under the
coordinate transformation%
\begin{equation}
S=\exp \left( \sigma _{1}x\right) ~,~\delta =\sigma _{2}\left( \rho x+\sqrt{%
1-\rho ^{2}}y\right)  \label{bs.09}
\end{equation}%
becomes%
\begin{equation}
F_{,xx}+F_{,yy}-\left( p_{1}x+p_{2}y+p_{3}\right) F_{,x}-\left(
q_{1}x+q_{2}y+q_{3}\right) F_{,y}-2F_{,t}=0  \label{bs.10}
\end{equation}%
in which the new parameters are expressed on the terms of the old ones according to
\begin{equation}
p_{1}=2\rho \frac{\sigma _{2}}{\sigma _{1}}~,~p_{2}=2\sqrt{1-\rho ^{2}}\frac{%
\sigma _{2}}{\sigma _{1}}~,~p_{3}=-2r,  \label{bs.11}
\end{equation}%
\begin{equation}
q_{1}=\frac{\kappa \sigma _{1}-\rho \sigma _{2}}{\sigma _{1}\sqrt{1-\rho ^{2}%
}},~q_{2}=\frac{\kappa \sigma _{1}-\rho \sigma _{2}}{\sigma _{1}}
\label{bs.12}
\end{equation}%
and%
\begin{equation}
q_{3}=-\frac{\left( \sigma _{1}^{2}\sigma _{2}\rho -2\sigma _{2}\rho
r+2\sigma _{1}\kappa \alpha -2\sigma _{1}\lambda \right) }{\sigma _{1}\sigma
_{2}\sqrt{1-\rho ^{2}}}.  \label{bs.13}
\end{equation}

The Lie symmetries for the autonomous two-factor model (\ref{bs.08}) have
been reported in \cite{Leach08}. However, for the convenience of the reader we
present the results.

\subsection{Lie symmetries of the autonomous equation}

Consider the infinitesimal one-parameter point transformation
\begin{eqnarray}
t^{\prime } &=&t+\varepsilon \xi ^{1}\left( t,x,y,F\right) ~,~x^{\prime
}=x+\varepsilon \xi ^{2}\left( t,x,y,F\right)  \label{bs.14} \\
y^{\prime } &=&y+\varepsilon \xi ^{3}\left( t,x,y,F\right) ~,~F^{\prime
}=y+\varepsilon \eta \left( t,x,y,F\right),  \label{bs.15}
\end{eqnarray}%
where $\varepsilon $ is an infinitesimal number so that $\varepsilon
^{2}\rightarrow 0$. From the transformation we define the generator $X$, as
\begin{equation}
X=\frac{\partial t^{\prime }}{\partial \varepsilon }\partial _{t}+\frac{%
\partial x^{\prime }}{\partial \varepsilon }\partial _{x}+\frac{\partial
y^{\prime }}{\partial \varepsilon }\partial _{y}+\frac{\partial F^{\prime }}{%
\partial \varepsilon }\partial _{F}  \label{bs.16}
\end{equation}%
or, equivalently,
\begin{equation}
X=\xi ^{1}\left( t,x,y,F\right) \partial _{t}+\xi ^{2}\left( t,x,y,F\right)
\partial _{x}+\xi ^{3}\left( t,x,y,F\right) \partial _{y}+\eta \left(
t,x,y,F\right) \partial _{F}.  \label{bs.17}
\end{equation}

The differential equation, $\Theta ,~$(\ref{bs.10}), is invariant under the
action of the one-parameter point transformation (\ref{bs.14})-(\ref{bs.15})
if there exists a function $\Lambda $ such that \cite{Olver,IbragB}%
\begin{equation}
X^{\left[ 2\right] }\Theta =\Lambda \Theta  \label{bs.18}
\end{equation}%
in which $X^{\left[ 2\right] }$ is the second prologation of $X$ defined in
the space $\left\{ t,x,y,F,F_{,x},F_{,y},F_{,xx},F_{,yy},F_{,xy}\right\} $.
When condition (\ref{bs.18}) holds, we say that $X$ is a Lie (point) symmetry
of $\Theta .$

Therefore from (\ref{bs.18}) we have the following Lie symmetries admitted
by equation (\ref{bs.10})
\begin{equation}
X_{t}=\partial _{t}~,~X_{F}=F\partial _{F}~,~X_{\infty }=f\left(
t,x,y\right) \partial _{f},
\end{equation}%
\begin{equation}
X_{1}=e^{c_{+}t}\left( a_{1}\partial _{x}+a_{2}\partial _{y}\right),
\end{equation}%
\begin{equation}
X_{2}=e^{c_{-}t}\left( a_{1}^{\prime }\partial _{x}+a_{2}^{\prime }\partial
_{y}\right),
\end{equation}%
\begin{equation}
X_{3}=e^{c_{+}t}\left( b_{1}\partial _{x}+b_{2}\partial _{y}+\left(
b_{3}x+b_{4}x+b_{5}\right) F\partial _{F}\right)
\end{equation}%
and
\begin{equation}
X_{4}=e^{c_{-}t}\left( b_{1}^{\prime }\partial _{x}+b_{2}^{\prime }\partial
_{y}+\left( b_{3}^{\prime }x+b_{4}^{\prime }x+b_{5}^{\prime }\right)
F\partial _{F}\right) .
\end{equation}%
The parameters $a_{1,2},~a_{1,2}^{\prime }$, $b_{1-5}$, $b_{1-5}^{\prime }$
and $c_{\pm }$ are functions of $p_{1-3}$~and $q_{1-3}$. $\ $The Lie
symmetries form the $\left\{ A_{1}\oplus _{s}W_{5}\right\} \oplus _{s}\infty
A_{1}$ Lie algebra. We note that for special cases of the parameters $%
p_{1-3},~q_{1-3}$, the representation of the admitted Lie symmetries of
equation (\ref{bs.10}) can be different. \ For instance, when all the
parameters $q_{1-3}\ $vanish, $q_{1-3}=0$, the Lie symmetries $X_{1-4}$
become%
\begin{equation}
X_{1}^{\prime }=p_{2}\partial _{x}-p_{1}\partial _{y}~,~X_{2}^{\prime }=e^{%
\frac{p_{1}}{2}t}\partial _{x},
\end{equation}%
\begin{equation}
X_{3}^{\prime }=\left( p_{1}p_{2}t+2p_{2}\right) \partial
_{x}-tp_{1}^{2}\partial _{y}+p_{1}^{2}yF\partial _{F}
\end{equation}%
and
\begin{equation}
X_{4}^{\prime }=e^{-\frac{p_{1}}{2}t}\left( \left(
p_{1}^{2}-p_{2}^{2}\right) \partial _{x}+2p_{1}p_{2}\partial
_{y}+p_{1}^{2}\left( p_{1}x+p_{2}y+p_{3}\right) F\partial _{F}\right).
\end{equation}%
For the remaining cases see \cite{Leach08}.

Below the nonautonomous two-factor model is defined and the group invariant
point transformations are derived.

\subsection{Lie symmetries of the nonautonomous equation}

We consider that the parameters $\sigma _{I},~\rho ,$ $r$, $\kappa ,~\alpha $
and $\lambda $ of (\ref{bs.08}) are well-defined functions of time. Without
loss of generality we can select a new time variable $\tau $ and eliminate,
for instance, the function $\sigma _{1}\left( t\right) $. Therefore we
select $\sigma _{1}=1.$

Under the time-depedent coordinate transformation, (\ref{bs.09}), the
two-factor model (\ref{bs.08}) has the following mathematical expression%
\begin{equation}
F_{,xx}+F_{,yy}-\left( P_{1}\left( t\right) x+P_{2}\left( t\right)
y+P_{3}\left( t\right) \right) F_{,x}-\left( Q_{1}\left( t\right)
x+Q_{2}\left( t\right) y+Q_{3}\left( t\right) \right) F_{,y}-2F_{,t}=0,
\label{bs.20a}
\end{equation}%
where now the new time-depedent parameters of the model are
\begin{equation}
P_{1}\left( t\right) =2\rho \sigma ~,~P_{2}\left( t\right) =2\sigma _{2}%
\sqrt{1-\rho ^{2}}~,~P_{3}\left( t\right) =1-2r\left( t\right),
\end{equation}%
\begin{equation}
Q_{1}\left( t\right) =-\frac{2\left( \rho \sigma _{2}\right) ^{2}+\left(
\rho \sigma _{2}\right) _{,t}+\rho \sigma _{2}\kappa }{\sigma _{2}\sqrt{%
1-\rho ^{2}}},
\end{equation}%
\begin{equation}
Q_{2}\left( t\right) =-\left( 2\rho \sigma _{2}+\kappa +2\frac{\sigma _{2,t}%
}{\sigma _{2}}\right) +\frac{2\rho _{2}\rho _{2,t}}{\sqrt{1-\rho ^{2}}}
\end{equation}%
and
\begin{equation}
Q_{3}\left( t\right) =-\left( \frac{\sigma _{2}\left( \rho -2r\rho \right)
-2\kappa \alpha +2\lambda }{\sigma _{2}\sqrt{1-\rho ^{2}}}\right) .
\end{equation}

Therefore, from the symmetry condition (\ref{bs.18}) for equation (\ref%
{bs.20a}), we find that the generic Lie symmetry vector is
\begin{eqnarray}
X_{G} &=&a\partial _{t}+\left( b_{1}+y\left( B_{2}+\frac{1}{4}aP_{2}-\frac{1%
}{4}aQ_{1}\right) +\frac{xa^{\prime }}{2}\right) \partial _{x}+  \nonumber \\
&&+\left( g-x\left( B_{2}+\frac{1}{4}aP_{2}-\frac{1}{4}aQ_{1}\right) +\frac{%
ya^{\prime }}{2}\right) \partial _{y}+  \nonumber \\
&&+\frac{1}{4}\left[ 4h+2xb_{1}P_{1}+2xgP_{2}+x^{2}(-P_{2}-Q_{1})\left(
B_{2}+\frac{1}{4}aP_{2}-\frac{1}{4}aQ_{1}\right) \right] F\partial _{F}+
\nonumber \\
&&+\frac{1}{4}\left[ 2x\left( B_{2}+\frac{1}{4}aP_{2}-\frac{1}{4}%
aQ_{1}\right) (yP_{1}-yQ_{2}-Q_{3})+x^{2}P_{1}a^{\prime }+2xyP_{2}a^{\prime }%
\right] F\partial _{F}+  \nonumber \\
&&+\frac{1}{4}\left[ xP_{3}a^{\prime }-4xb_{1}^{\prime 2}aP_{1}^{\prime
}+2xyaP_{2}^{\prime }+2xaP_{3}^{\prime }\right] F\partial _{F}+  \nonumber \\
&&+\frac{1}{4}\left[ -4xy\left( \frac{1}{4}P_{2}a^{\prime }-\frac{1}{4}%
Q_{1}a^{\prime }+\frac{1}{4}aP_{2}^{\prime }-\frac{1}{4}aQ_{1}^{\prime
}\right) -x^{2}a^{\prime \prime }\right] F\partial _{F}+  \nonumber \\
&&+\frac{1}{4}\left[ 2yb_{1}Q_{1}+2yP_{3}\left( B_{2}+\frac{1}{4}aP_{2}-%
\frac{1}{4}aQ_{1}\right) +y^{2}(P_{2}+Q_{1})\left( B_{2}+\frac{1}{4}aP_{2}-%
\frac{1}{4}aQ_{1}\right) \right] F\partial _{F}+  \nonumber \\
&&+\frac{1}{4}\left[ 2ygQ_{2}+y^{2}Q_{2}a^{\prime }+yQ_{3}a^{\prime
}-4yg^{\prime 2}aQ_{2}^{\prime }+2yaQ_{3}^{\prime 2}a^{\prime \prime }\right]
F\partial _{F},  \label{bs.26bb}
\end{eqnarray}
where $B_{2}$ is constant, $a=a\left( t\right) ,~b_{1}=b_{1}\left( t\right)
,~f=f\left( t\right) $ and $g=g\left( t\right) $, given by the system of
equations of Appendix \ref{ap1}. Furthermore, from the generic vector field (%
\ref{bs.26bb}) and the system of Appendix \ref{ap1}, we have that the
nonautonomous two-factor model of commodities is invariant under the $%
\left\{ A_{1}\oplus _{s}W_{5}\right\} \oplus _{s}\infty A_{1}$ Lie algebra,
the same algebra as the autonomous model but in a different representation.

We continue our analysis with the two-dimensional Black-Scholes equation.

\section{The two-dimensional Black-Scholes equation}

\label{2DBS}

Consider a basket containing two assets the prices of which are $S_{1}$ and $%
S_{2}$ and that the the prices of the underlying assets obey the system of
stochastic differential equations,
\begin{equation}
dS_{I,t}=S_{I,t}\left( \mu _{I}dt+\frac{\sigma _{I}}{\sqrt{1+\rho ^{2}}}%
\left( dW_{I,t}+\rho dW_{J,t}\right) \right),  \label{bs.20}
\end{equation}%
where $I,J=1,2$, $\ I\neq J$, and $W_{I,t}$ are two independent standard
Brownian motions. Let $u=u\left( t,S_{1},S_{2}\right) $ be the payoff
function on a European option on this two-asset basket. Then the evolution
equation which $u$ satisfies is an $(1+2)$ linear evolution equation given
by \cite{finBook}%
\begin{equation}
\frac{1}{2}\sigma _{1}^{2}u_{,11}+\rho \sigma _{1}\sigma _{2}u_{,12}+\frac{1%
}{2}\sigma _{2}^{2}u_{,22}-rS_{1}u_{,1}-rS_{2}u_{,2}-ru+u_{,t}=0
\label{bs.21}
\end{equation}%
with the terminal condition $u\left( T,S_{1},S_{2}\right) =U$, when $t=T.$

Equation (\ref{bs.21}) is a generalisation of the BS equation and it is
called the two-dimensional BS equation. \ The Lie symmetry analysis of (\ref%
{bs.01}) has been presented in \cite{ibrag} and recently a Lie symmetry
analysis for equation (\ref{bs.01}), with a general potential function, was
performed in \cite{dimas2}. The algebraic properties of the autonomous
form of (\ref{bs.21}) have been studied in \cite{consta} and it was found
that equation (\ref{bs.21}) is invariant under a seven-dimensional Lie
algebra, plus the infinite number of solution symmetries. As we see below, the
analysis of the autonomous equation (\ref{bs.21}) in \cite{consta} is not
complete. In particular we find that it is maximally symmetric, \textit{ie}
invariant under a nine-dimensional Lie algebra, plus the infinite number of solution
symmetries. In \cite{consta} the authors considered the following equation%
\begin{equation}
\frac{1}{2}\sigma _{1}^{2}u_{,11}+\rho \sigma _{1}\sigma _{2}u_{,12}+\frac{1%
}{2}\sigma _{2}^{2}u_{,22}-\mu _{1}S_{1}u_{,1}-\mu
_{2}S_{2}u_{,2}-ku+u_{,t}=0  \label{bs.22}
\end{equation}%
which reduces to (\ref{bs.21}) when $\mu _{1}=\mu _{2}=k=r.$

Below we determine the Lie symmetries of equation (\ref{bs.22}) for the
autonomous and nonautonomous system.

\subsection{Lie symmetries of the autonomous equation}

We introduce the coordinate transformation%
\begin{equation}
S_{1}=\exp \left( \sigma _{1}x\right) ~,~S_{2}=\exp \left( \sigma _{2}\rho
x+\sigma _{2}\sqrt{1-\rho ^{2}}y\right)  \label{bs.23}
\end{equation}%
under which equation (\ref{bs.22}) becomes%
\begin{equation}
u_{,xx}+u_{,yy}-\phi _{1}u_{,x}-\phi _{2}u_{,y}-2ku+2u_{,t}=0,  \label{bs.24}
\end{equation}%
where now the new constants, $\phi _{1}$ and $\phi _{2}$, are
\begin{equation}
\phi _{I}=\frac{\sigma _{1}^{2}+2\mu _{I}}{\sigma _{I}}.  \label{bs.25}
\end{equation}

On application of the Lie symmetry condition (\ref{bs.18}) for (%
\ref{bs.24}) we find that the Lie symmetry vectors are
\begin{equation}
X_{t}=\partial _{t}~,~X_{u}=F\partial _{u}~,~X_{\infty }=f\left(
t,x,y\right) \partial _{u},
\end{equation}%
\[
X_{1}=\partial _{x}~,~X_{2}=t\partial _{x}+\frac{1}{2}x`x\left( x+\phi
_{1}t\right) u\partial _{u},
\]%
\begin{equation}
X_{3}=\partial _{y}~,~X_{4}=t\partial _{y}+\frac{1}{2}\left( y+\phi
_{2}t\right) u\partial _{u},
\end{equation}%
\begin{equation}
X_{5}=y\partial _{x}-x\partial _{y}+\frac{1}{2}\left( \phi _{1}y-\phi
_{2}x\right) u\partial _{u},
\end{equation}%
\begin{equation}
X_{6}=2t\partial _{t}+x\partial _{x}+y\partial _{y}+\frac{1}{2}\left( \phi
_{1}x+\phi _{2}y+t\left( \phi _{1}^{2}+\phi _{2}^{2}+8k\right) \right)
u\partial _{u}
\end{equation}%
and
\begin{equation}
X_{7}=t^{2}\partial _{t}+tx\partial _{x}+ty\partial _{y}+\frac{1}{4}\left(
x^{2}+y^{2}+t^{2}\left( \phi _{1}^{2}+\phi _{2}^{2}+8k\right) +2t\left( \phi
_{1}x+\phi _{2}y-2\right) \right) u\partial _{u}.
\end{equation}%
which are $8+1+\infty $ symmetries.  This is the admitted group invariant
algebra of the two-dimensional Heat Equation, that is, $\left\{ \left\{
sl\left( 2,R\right) \oplus _{s}so\left( 2\right) \right\} \oplus
_{s}W_{5}\right\} \oplus _{s}\infty A_{1}$. Hence the two-dimensional BS
equation (\ref{bs.22}) is maximally symmetric and equivalent with the
two-dimensional Heat and Schr\"{o}dinger equations \cite{AnJGP}. This result
does not hold for the two-factor model of commodities. An analysis
does hold when in (\ref{bs.22}), $\mu _{1}=\mu _{2}=k=r$; that is, for
equation (\ref{bs.21}).

When we apply the transformations
\begin{equation}
t=-\frac{1}{2}T~,~x=\bar{x}-\frac{1}{2}\phi _{1}t~
\end{equation}%
and
\begin{equation}
~\bar{y}=y-\frac{1}{2}\phi _{2}t~,~u=e^{2kt}v\left( t,x,y\right)
\end{equation}%
to (\ref{bs.24}), the equation becomes
\begin{equation}
v_{,\bar{x}\bar{x}}+v_{,\bar{y}\bar{y}}-v_{,t}=0  \label{h2d}
\end{equation}%
which is the two-dimensional Heat conduction equation.

We proceed to the determination of the Lie symmetries for the
nonautonomous equation (\ref{bs.22}).

\subsection{Lie symmetries of the nonautonomous equation}

We take the parameters, $\sigma _{I},~\rho ,$ $\mu _{I}$ and $k, $ of (\ref%
{bs.22}) to be well-defined functions of time. Moreover without loss of
generality we select $\sigma _{1}\left( t\right) =1.$

We apply the time-dependent transformation (\ref{bs.23}) to (\ref{bs.22})
and we have%
\begin{equation}
u_{,xx}+u_{,yy}-P_{1}\left( t\right) u_{,x}-\left( Q_{1}\left( t\right)
x+Q_{2}\left( t\right) y+Q_{3}\left( t\right) \right) u_{,y}-2k\left(
t\right) u+2u_{,t}=0  \label{bs.26}
\end{equation}%
in which%
\begin{equation}
P_{1}\left( t\right) =1+2\mu _{1}\left( t\right) ~,~Q_{1}\left( t\right) =%
\frac{2\left( \rho \sigma _{2}\right) _{,t}}{\sigma _{2}\sqrt{1-\rho ^{2}}},
\end{equation}%
\begin{equation}
Q_{2}\left( t\right) =-\frac{2\left( \sigma _{2,t}\rho ^{2}+\sigma _{2}\rho
\rho _{,t}-\sigma _{2,t}\right) }{\sigma _{2}\left( 1-\rho ^{2}\right) }
\end{equation}%
and
\begin{equation}
Q_{3}\left( t\right) =\frac{\sigma _{2}\left( \sigma _{2}-\rho -2\mu
_{2}\rho \right) +2\mu _{2}}{\sigma _{2}\sqrt{1-\rho ^{2}}}
\end{equation}

From the symmetry condition (\ref{bs.18}) for equation (\ref{bs.26}%
) we find that the generic Lie symmetry vector has the following
mathematical expression
\begin{eqnarray}
X_{G} &=&a\partial _{t}+\left( b_{1}+y\left( B_{2}+\frac{1}{4}aQ_{1}\right) +%
\frac{xa^{\prime }}{2}\right) \partial _{x}+\left( f-x\left( B_{2}+\frac{1}{4%
}aQ_{1}\right) +\frac{ya^{\prime }}{2}\right) \partial _{y}+  \nonumber \\
&&+\frac{1}{4}\left[ 4g+\left( -x^{2}Q_{1}\left( B_{2}+\frac{1}{4}%
aQ_{1}\right) -2x\left( B_{2}+\frac{1}{4}aQ_{1}\right) (yQ_{2}+Q_{3})\right) %
\right] u\partial _{u}+  \nonumber \\
&&+\frac{1}{4}\left[ xP_{1}a^{\prime }+4xb_{1}^{\prime }+2xaP_{1}^{\prime
}+x^{2}a^{\prime \prime }+4xy\left( \frac{1}{4}Q_{1}a^{\prime }+\frac{1}{4}%
aQ_{1}^{\prime }\right) \right] u\partial _{u}+  \nonumber \\
&&+\frac{1}{4}\left[ +2yb_{1}Q_{1}+2yP_{1}\left( B_{2}+\frac{1}{4}%
aQ_{1}\right) +y^{2}Q_{1}\left( B_{2}+\frac{1}{4}aQ_{1}\right) \right]
u\partial _{u}+  \nonumber \\
&&+\frac{1}{4}\left[ 2yfQ_{2}+y^{2}Q_{2}a^{\prime }+yQ_{3}a^{\prime
}+4yf^{\prime 2}aQ_{2}^{\prime }+2yaQ_{3}^{\prime 2}a^{\prime \prime }\right]
u\partial _{u},  \label{bs.26d}
\end{eqnarray}%
\newline
where $B_{2}$ is a constant, $a=a\left( t\right) ,~b_{1}=b_{1}\left(
t\right) ,~f=f\left( t\right) $ and $g=g\left( t\right) $ which given by the
system of differential equations of Appendix \ref{ap2}. Furthermore, from (%
\ref{bs.26d}) and the system of Appendix \ref{ap2}, we observe that the
nonautonomous equation (\ref{bs.21}) is invariant under the group of
transformations in which the generators form the $\left\{ \left\{
A_{1}\oplus _{s}so\left( 2\right) \right\} \oplus _{s}W_{5}\right\} \oplus
_{s}\infty A_{1}$ Lie algebra. \ Below we consider a special case for which $%
\sigma _{1}\left( t\right) \simeq \sigma _{2}\left( t\right) $ and $\rho
=const.$

\subsubsection{Special Case: $\protect\rho =const$ and $\protect\sigma %
_{1}\left( t\right) \simeq \protect\sigma _{2}\left( t\right) $}

As a special case of the nonautonomous equation (\ref{bs.22}) we consider $%
\sigma _{2}\left( t\right) =\sigma _{0}\sigma _{1}\left( t\right) $, where $%
\sigma _{0}$ is a constant and $\rho \left( t\right) $ is a constant.
The nonautonomous two-dimensional BS equation becomes%
\begin{equation}
\sigma _{1}^{2}\left( t\right) \left( \frac{1}{2}u_{,11}+\rho \sigma
_{0}u_{,12}+\frac{1}{2}\sigma _{0}^{2}u_{,22}\right) -\mu _{1}\left(
t\right) S_{1}u_{,1}-\mu _{2}\left( t\right) S_{2}u_{,2}-k\left( t\right)
u+u_{,t}=0,  \label{bs.27}
\end{equation}%
where without loss of generality we can select $\sigma _{1}\left( t\right)
=1 $. Under the transformation (\ref{bs.23}) equation (\ref{bs.27}) becomes%
\begin{equation}
u_{,xx}+u_{,yy}-\Lambda _{1}\left( t\right) u_{,x}-\Lambda _{2}\left(
t\right) u_{,y}-2k\left( t\right) u+2u_{,t}=0,  \label{bs.28}
\end{equation}%
where the new functions $\Lambda _{1}\left( t\right) ,~\Lambda _{2}\left(
t\right) $ are defined as%
\begin{equation}
\Lambda _{1}\left( t\right) =1+2\mu _{1}\left( t\right)
\end{equation}%
and
\begin{equation}
\Lambda _{2}\left( t\right) =\frac{\sigma _{0}\left( \sigma _{0}-\rho -2\mu
_{2}\left( t\right) \rho \right) +2\mu _{2}\left( t\right) }{\sigma _{0}%
\sqrt{1-\rho ^{2}}}.
\end{equation}

From the symmetry condition (\ref{bs.18}) for equation (\ref{bs.26}) the
following symmetry vectors arise%
\begin{equation}
X_{u}=u\partial _{u}~~,~X_{\infty }=f\left( t,x,y\right) \partial _{F},
\end{equation}%
\begin{equation}
Z_{1}=\partial _{x}~~,~~Z_{2}=t\partial _{x}+\left( \frac{1}{2}\int \Lambda
_{1}dt+x\right) u\partial _{u},
\end{equation}%
\begin{equation}
Z_{3}=\partial _{y}~\ ,~Z_{4}=t\partial _{y}+\left( \frac{1}{2}\int \Lambda
_{2}dt+y\right) u\partial _{u},
\end{equation}%
\begin{equation}
Z_{5}=\left( y+\frac{1}{2}\int \Lambda _{2}dt\right) \partial _{x}-\left( x+%
\frac{1}{2}\int \Lambda _{1}dt\right) \partial _{y}+\frac{1}{2}\left(
\Lambda _{1}y-\frac{1}{2}\Lambda _{2}x\right) u\partial _{u},
\end{equation}%
\begin{equation}
Z_{6}=\partial _{t}-\frac{1}{2}\Lambda _{1}\partial _{x}-\frac{1}{2}\Lambda
_{2}\partial _{y}+ku\partial _{u},
\end{equation}%
\begin{equation}
Z_{7}=2t\partial _{t}+\left( x-\frac{1}{2}\int \Lambda _{1}dt-\int t\Lambda
_{1}dt\right) \partial _{x}+\left( y-\frac{1}{2}\int \Lambda _{2}dt-\int
t\Lambda _{2}dt\right) \partial _{y}+tku\partial _{u}
\end{equation}%
and
\begin{eqnarray}
Z_{8} &=&t^{2}\partial _{t}+\left( tx-\frac{1}{2}\int \int \left(
t^{2}\Lambda _{1,tt}+3t\Lambda _{1,t}\right) dt\right) \partial _{x}+\left(
ty-\frac{1}{2}\int \int t^{2}\Lambda _{2,tt}+3t\Lambda _{2,t}\right)
\partial _{y}+  \nonumber \\
&&+\left[ -\frac{1}{2}x\left( \int t^{2}\Lambda _{1,tt}dt+3\int t\Lambda
_{1,t}dt-t^{2}\Lambda _{1,t}-t\Lambda _{1}-x\right) -\right] u\partial _{u}+
\nonumber \\
&&+\left[ -\frac{1}{2}y\left( \int t^{2}\Lambda _{2,tt}dt+3\int t\Lambda
_{2,t}dt-t^{2}\Lambda _{2,t}-t\Lambda _{2}-y\right) \right] u\partial _{u}+
\nonumber \\
&&+\frac{1}{4}\left[ 4t\left( t-1\right) -\int \Lambda _{1}\left( \int
t^{2}\Lambda _{1,tt}dt\right) dt-\int \Lambda _{2}\left( \int t^{2}\Lambda
_{2,tt}dt\right) dt\right] u\partial _{u}+  \nonumber \\
&&+\frac{1}{4}\left[ -3\int \Lambda _{1}\int t\Lambda _{1,t}dt-3\int \Lambda
_{2}\int t\Lambda _{2,t}dt\right] u\partial _{u}  \nonumber \\
&&+\frac{1}{4}\left[ \int t^{2}\Lambda _{1}\Lambda _{1,t}+\int t^{2}\Lambda
_{2}\Lambda _{2,t}+\int t\left( \Lambda _{1}^{2}+\Lambda _{2}^{2}\right) dt%
\right] u\partial _{u}.
\end{eqnarray}

Hence the nonautonomous equation (\ref{bs.27}) is maximally symmetric, just as
the autonomous two-dimensional BS equation,\ in contrast to the
nonautonomous equation (\ref{bs.26}) which is invariant under another group
of point transformations.

Moreover equation (\ref{bs.28}) can be written in the form of (\ref{h2d})
and the transformation which does that is%
\begin{equation}
t=-\frac{1}{2}T~,~u=e^{2kt}v\left( t,x,y\right) ,
\end{equation}%
and
\begin{equation}
x=\bar{x}-\frac{1}{2}\int \Lambda _{1}dt~,~y=\bar{y}-\frac{1}{2}\int \Lambda
_{2}dt~.
\end{equation}

Below we discuss our results and draw our conclusions

\section{Conclusions}

\label{con}

The purpose of this work is to study the algebraic properties of nonautonomous $(1+2)~$evolution equations in financial mathematics.
Specifically we examined the relation among the admitted group of invariant
transformations between the autonomous and the nonautonomous equations of
the two-factor model of commodities and of the two-dimensional BS equation
was performed.

For the two-factor model of commodities we proved that the autonomous and
the nonautonomous equations are invariant under the same group of point
transformations in which the generators form the $\left\{ A_{1}\oplus
_{s}W_{5}\right\} \oplus _{s}\infty A_{1}$ Lie algebra.

As far as the autonomous two-dimensional BS equation is concerned, we proved
that it is maximally symmetric and admits as Lie symmetries the generators of
the Lie algebra \ $\left\{ \left\{ sl\left( 2,R\right) \oplus _{s}so\left(
2\right) \right\} \oplus _{s}W_{5}\right\} \oplus _{s}\infty A_{1}$ This
corrects the existing result in the literature. \ However, the admitted Lie
symmetries of the nonautonomous two-dimensional BS equation form a different
Lie algebra than that of the autonomous equation and is of lower
dimension. Specifically the admitted Lie algebra is $\left\{ \left\{
A_{1}\oplus _{s}so\left( 2\right) \right\} \oplus _{s}W_{5}\right\} \oplus
_{s}\infty A_{1}$. That result differs from that for the model of
commodities for which the autonomous and the nonautonomous equations are
invariant under the same group of transformations, namely $\left\{
A_{1}\oplus _{s}W_{5}\right\} \oplus _{s}\infty A_{1}.$

In the case for which $\rho =const$ and $\sigma _{1}\left( t\right) \simeq
\sigma _{2}\left( t\right) $, the two-dimensional BS equation is maximally
symmetric. In order to understand why we have this special case consider the
general $\left( 1+n\right) $ evolution equation ( We use the Einstein
summation convention).
\begin{equation}
A^{ij}\left( t,x^{k}\right) u_{ij}+B^{i}\left( t,x^{k}\right) u_{,i}+f\left(
t,x^{k},u\right) =u_{,t}.  \label{con.1}
\end{equation}

If $X=\xi ^{t}\partial _{t}+\xi ^{i}\partial _{i}+\eta \partial _{u}$ is the
generator of a Lie symmetry vector, one of the symmetry conditions can be
written as%
\begin{equation}
\mathcal{L}_{\xi ^{\alpha }}A^{ij}=-2\psi A^{ij},  \label{con.2}
\end{equation}%
where $\psi $ is a function of $t$ only, and $\alpha =1,2,...,n,t$.
Therefore from (\ref{con.2}) we have that
\begin{equation}
\mathcal{L}_{\xi ^{i}}A^{ij}=-2\psi A^{ij}-A_{,t}^{ij}\xi ^{t}.
\label{con.3}
\end{equation}

From (\ref{con.3}) we have that, when $A_{,t}^{ij}=0$, the Lie symmetries of
(\ref{con.1}) are generated by the Homothetic Algebra of $A_{ij}$. However,
that is not true when $A_{,t}^{ij}\neq 0$ and new possible generators arise.
\ In the $\left( 1+1\right) $ equations, \textit{ie} (\ref{bs.01})~and (\ref%
{bs.02}), when $\sigma =\sigma \left( t\right) $, as we discussed above, we
can always perform a time (coordinate) transformation and cause the second
derivatives to be time-independent. Therefore, in order to apply this method
to the two-dimensional systems, we have to select $\rho =const$ and $\sigma
_{1}\left( t\right) \simeq \sigma _{2}\left( t\right) $ so that at the end
the components of the second derivatives can be seen as time-independent.

Furthermore we remark that we performed a reduction on the two nonautonomous
equations (\ref{bs.08}) and (\ref{bs.21}) by using the Lie symmetries (\ref%
{bs.26bb}) and (\ref{bs.26d}), respectively, for $a\left( t\right) =0$. We
found that the reduced equations, which are $\left( 1+1\right) $ evolution
equations, are maximally symmetric. This is the same result as is to be
found in the case of the autonomous two-factor model \cite{Leach08}.

As a final application consider the nonautonomous two-dimensional BS
equation (\ref{bs.28}). From the application of the invariant functions of
the Lie symmetries $\left\{ Z_{1}+c_{1}X_{u},~Z_{3}+c_{2}X_{u}\right\} $ we
have the solution $u\left( t,x,y\right) =w\left( t\right) \exp \left(
c_{1}x+c_{2}y\right) $, where%
\begin{equation}
w\left( t\right) =\exp \left( \frac{1}{2}\int \left( 2k\left( t\right)
-\left( c_{1}^{2}+c_{2}^{2}\right) +\Lambda _{1}\left( t\right)
c_{1}+\Lambda _{2}\left( t\right) c_{2}\right) dt\right) .  \label{con.5}
\end{equation}

\begin{figure}[tbp]
\centering\includegraphics[height=6cm]{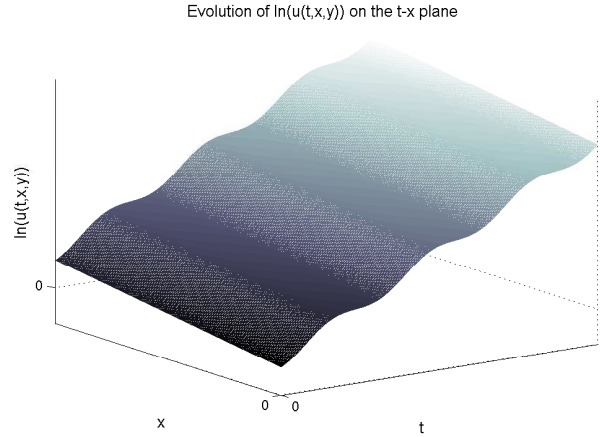}
\caption{Qualitative evolution of the solution $u\left( t,x,y\right) $ for
the nonautonomous two-dimensional Black-Scholes-Merton equation (\protect\ref%
{bs.21}) in the \textquotedblleft t-x\textquotedblright\ plane, when $%
\protect\sigma _{1},\protect\sigma _{2},\protect\rho $ are constants and $%
r\left( t\right) =r_{0}+\protect\varepsilon \sin \left( \protect\omega %
t\right) $. }
\label{fig3}
\end{figure}

In the case for which $\mu _{1}\left( t\right) =\mu _{2}\left( t\right)
=k\left( t\right) =r\left( t\right) $ and $r\left( t\right)
=r_{0}+\varepsilon \sin \left( \omega t\right) $, $\omega $ ,$\varepsilon$ and $
r_{0}~$are constants, the solution of the nonautonomous two-dimensional BS
equation for the \textquotedblleft $t-x$\textquotedblright\ plane is given
in figure \ref{fig3}. We observe that in the $t-$direction, function $%
u\left( t,x,y\right) $ has periodic behavior along the line $f\left(
t\right) \simeq t$ with period $\omega $.

{The implication of the results of the present analysis is that for
the two-factor model of commodities, the autonomous and the nonautonomous
problem share the same static solutions, that is, the differences follow
only from the time-dependent terms of the solution. However, that is not true
for the two-dimensional Black-Scholes Equation in which the nonautonomous
equation in general is not maximally symmetric and does not share the same
number of static solutions with that of the autonomous equation. On the
other hand we found that if and only if the time-dependence of the two
volatilities }$\sigma _{1}\left( t\right) ,~\sigma _{2}\left( t\right) $%
\textbf{\ are the same, i.e., }$\frac{\sigma _{1}\left( t\right) }{\sigma
_{2}\left( t\right) }=const${, and that the correlation factor }$\rho
${\ is constant then the nonautonomous Black-Scholes shares the same
static solutions, i.e. static evolution, with the autonomous equation. }

{The results of this analysis are important in the sense that by
starting from the autonomous equation and with the use of coordinate
transformations and only someone can analyse models with time-varying
constants. On the other hand starting from real data and with the use of
coordinate transformations to see if the data are well described from the
autonomous system, and vice verca. The situation is not different from that
which one finds on the relation between the free particle and harmonic
oscillator. In order to demonstrate that, if we plot the time-position
diagram of the mathematical pendulum, where we measure the distance and the
time with nonlinear instruments, the graph will be a straight line, which
describes the motion of the free particle. }

In a forthcoming work we intend to extend our analysis to the case where the
free parameters of the models are space-dependent. Such an analysis it is in
progress and will be published elsewhere.

{\textbf{Acknowledgments:} The research of AP was supported by FONDECYT grant no. 3160121. RMM thanks the National Research Foundation of the Republic of
South Africa for the granting of a postdoctoral fellowship with grant number
93183 while this work was being undertaken}
\appendix 


\section*{\noindent Nonautonomous two-factor model of commodities}

\vspace{6pt}

\label{ap1}

In this Appendix we give the differential equations which the functions$%
~a\left( t\right) ,~b_{1}\left( t\right) ,~h\left( t\right) $ and $g\left(
t\right) $ of the generic symmetry vector (\ref{bs.26bb}) of the
nonautonomous two-factor model of commodities satisfy. For the derivation of
the system the symbolic package SYM of Mathematica has been used \cite%
{Dimas05a,Dimas06a,Dimas08a}.

The system is:

\begin{eqnarray}
0 &=&-\frac{1}{2}b_{1}P_{1}P_{3}-\frac{1}{2}gP_{2}P_{3}-\frac{1}{2}%
b_{1}Q_{1}Q_{3}-\frac{1}{2}gQ_{2}Q_{3}+  \nonumber \\
&&+\frac{1}{2}P_{1}a^{\prime }-\frac{1}{4}P_{3}^{2}a^{\prime }+\frac{1}{2}%
Q_{2}a^{\prime }-\frac{1}{4}Q_{3}^{2}a^{\prime }+P_{3}b_{1}^{\prime }+
\nonumber \\
&&+Q_{3}g^{\prime }-2h^{\prime }+\frac{1}{2}aP_{1}^{\prime }-\frac{1}{2}%
aP_{3}P_{3}^{\prime }+\frac{1}{2}aQ_{2}^{\prime }-\frac{1}{2}%
aQ_{3}Q_{3}^{\prime }-a^{\prime \prime },
\end{eqnarray}

\begin{eqnarray}
0 &=&-\frac{1}{2}b_{1}P_{1}^{2}-\frac{1}{2}gP_{1}P_{2}+\frac{1}{2}%
B_{2}P_{2}P_{3}+\frac{1}{8}aP_{2}^{2}P_{3}-\frac{1}{8}aP_{2}P_{3}Q_{1}+
\nonumber \\
&&-\frac{1}{2}b_{1}Q_{1}^{2}-\frac{1}{2}gQ_{1}Q_{2}+\frac{1}{2}%
B_{2}Q_{2}Q_{3}+\frac{1}{8}aP_{2}Q_{2}Q_{3}-\frac{1}{8}aQ_{1}Q_{2}Q_{3}+
\nonumber \\
&&-\frac{3}{4}P_{1}P_{3}a^{\prime }-\frac{3}{4}Q_{1}Q_{3}a^{\prime
}-P_{2}g^{\prime }+Q_{1}g^{\prime }-b_{1}P_{1}^{\prime }-\frac{1}{2}%
aP_{3}P_{1}^{\prime }+  \nonumber \\
&&-gP_{2}^{\prime }-\frac{1}{2}aP_{1}P_{3}^{\prime }-\frac{3a^{\prime
}P_{3}^{\prime }}{2}-\frac{1}{2}aQ_{3}Q_{1}^{\prime }+  \nonumber \\
&&+B_{2}Q_{3}^{\prime }+\frac{1}{4}aP_{2}Q_{3}^{\prime }-\frac{3}{4}%
aQ_{1}Q_{3}^{\prime }+2b_{1}^{\prime \prime }-aP_{3}^{\prime \prime },
\end{eqnarray}

\begin{eqnarray}
0 &=&-\frac{1}{2}b_{1}P_{1}P_{2}-\frac{1}{2}gP_{2}^{2}-\frac{1}{2}%
B_{2}P_{1}P_{3}-\frac{1}{8}aP_{1}P_{2}P_{3}+\frac{1}{8}aP_{1}P_{3}Q_{1}+
\nonumber \\
&&-\frac{1}{2}b_{1}Q_{1}Q_{2}-\frac{1}{2}gQ_{2}^{2}-\frac{1}{2}%
B_{2}Q_{1}Q_{3}-\frac{1}{8}aP_{2}Q_{1}Q_{3}+\frac{1}{8}aQ_{1}^{2}Q_{3}+
\nonumber \\
&&-\frac{3}{4}P_{2}P_{3}a^{\prime }-\frac{3}{4}Q_{2}Q_{3}a^{\prime
}+P_{2}b_{1}^{\prime }-Q_{1}b_{1}^{\prime }-\frac{1}{2}aP_{3}P_{2}^{\prime }+
\nonumber \\
&&-B_{2}P_{3}^{\prime }-\frac{3}{4}aP_{2}P_{3}^{\prime }+\frac{1}{4}%
aQ_{1}P_{3}^{\prime }-b_{1}Q_{1}^{\prime }-gQ_{2}^{\prime }+  \nonumber \\
&&-\frac{1}{2}aQ_{3}Q_{2}^{\prime }-\frac{1}{2}aQ_{2}Q_{3}^{\prime }-\frac{%
3a^{\prime }Q_{3}^{\prime }}{2}+2g^{\prime \prime }-aQ_{3}^{\prime \prime }
\end{eqnarray}%
and
\begin{eqnarray}
0 &=&B_{2}P_{1}P_{2}+\frac{1}{4}aP_{1}P_{2}^{2}-\frac{1}{4}%
aP_{1}P_{2}Q_{1}+B_{2}Q_{1}Q_{2}+  \nonumber \\
&&+\frac{1}{4}aP_{2}Q_{1}Q_{2}-\frac{1}{4}aQ_{1}^{2}Q_{2}-\frac{1}{2}%
P_{1}^{2}a^{\prime }+\frac{1}{2}P_{2}^{2}a^{\prime }-\frac{1}{2}%
Q_{1}^{2}a^{\prime }+  \nonumber \\
&&+\frac{1}{2}Q_{2}^{2}a^{\prime }-\frac{1}{2}aP_{1}P_{1}^{\prime
}-a^{\prime }P_{1}^{\prime }+B_{2}P_{2}^{\prime }+\frac{3}{4}%
aP_{2}P_{2}^{\prime }-\frac{1}{4}aQ_{1}P_{2}^{\prime }+  \nonumber \\
&&+B_{2}Q_{1}^{\prime }+\frac{1}{4}aP_{2}Q_{1}^{\prime }-\frac{3}{4}%
aQ_{1}Q_{1}^{\prime }+\frac{1}{2}aQ_{2}Q_{2}^{\prime }+a^{\prime
}Q_{2}^{\prime }-\frac{1}{2}aP_{1}^{\prime \prime }+\frac{1}{2}%
aQ_{2}^{\prime \prime }.
\end{eqnarray}

\section{Nonautonomous two-dimensional Black-Scholes}

\label{ap2}

In this Appendix we give the differential equations which the functions$%
~a\left( t\right) ,~b_{1}\left( t\right) ,~f\left( t\right) $ and $g\left(
t\right) $ of the generic symmetry vector (\ref{bs.26d}) of the
nonautonomous two-dimensional Black-Scholes Equation satisfy.

The system is:

\begin{eqnarray}
0 &=&-\frac{1}{2}b_{1}Q_{1}Q_{3}-\frac{1}{2}fQ_{2}Q_{3}-2ka^{\prime }-\frac{1%
}{4}P_{1}^{2}a^{\prime }+  \nonumber \\
&&+\frac{1}{2}Q_{2}a^{\prime }-\frac{1}{4}Q_{3}^{2}a^{\prime
}-P_{1}b_{1}^{\prime }-Q_{3}f^{\prime }+2g^{\prime }-2ak^{\prime }+
\nonumber \\
&&-\frac{1}{2}aP_{1}P_{1}^{\prime }+\frac{1}{2}aQ_{2}^{\prime }-\frac{1}{2}%
aQ_{3}Q_{3}^{\prime }+a^{\prime \prime },
\end{eqnarray}

\begin{eqnarray}
0 &=&-\frac{1}{2}b_{1}Q_{1}^{2}-\frac{1}{2}fQ_{1}Q_{2}+\frac{1}{2}%
B_{2}Q_{2}Q_{3}+\frac{1}{8}aQ_{1}Q_{2}Q_{3}-\frac{3}{4}Q_{1}Q_{3}a^{\prime }+
\nonumber  \label{1.3} \\
&&-Q_{1}f^{\prime }+\frac{3a^{\prime }P_{1}^{\prime }}{2}-\frac{1}{2}%
aQ_{3}Q_{1}^{\prime }-B_{2}Q_{3}^{\prime }-\frac{3}{4}aQ_{1}Q_{3}^{\prime
}+2b_{1}^{\prime \prime }+aP_{1}^{\prime \prime },
\end{eqnarray}

\begin{eqnarray}
0 &=&-\frac{1}{2}b_{1}Q_{1}Q_{2}-\frac{1}{2}fQ_{2}^{2}-\frac{1}{2}%
B_{2}Q_{1}Q_{3}-\frac{1}{8}aQ_{1}^{2}Q_{3}-\frac{3}{4}Q_{2}Q_{3}a^{\prime }+
\nonumber  \label{1.5} \\
&&+Q_{1}b_{1}^{\prime }+B_{2}P_{1}^{\prime }+\frac{1}{4}aQ_{1}P_{1}^{\prime
}+b_{1}Q_{1}^{\prime }+fQ_{2}^{\prime }+  \nonumber \\
&&-\frac{1}{2}aQ_{3}Q_{2}^{\prime }-\frac{1}{2}aQ_{2}Q_{3}^{\prime }+\frac{%
3a^{\prime }Q_{3}^{\prime }}{2}+2f^{\prime \prime }+aQ_{3}^{\prime \prime },
\end{eqnarray}

\begin{eqnarray}
0 &=&-\frac{1}{2}B_{2}Q_{1}^{2}-\frac{1}{8}aQ_{1}^{3}+\frac{1}{2}%
B_{2}Q_{2}^{2}+\frac{1}{8}aQ_{1}Q_{2}^{2}+  \nonumber  \label{1.6} \\
&&-Q_{1}Q_{2}a^{\prime }-\frac{1}{2}aQ_{2}Q_{1}^{\prime }+a^{\prime
}Q_{1}^{\prime }-B_{2}Q_{2}^{\prime }-\frac{3}{4}aQ_{1}Q_{2}^{\prime }+\frac{%
1}{2}aQ_{1}^{\prime \prime }
\end{eqnarray}%
and

\begin{eqnarray}
0 &=&-B_{2}Q_{1}Q_{2}-\frac{1}{4}aQ_{1}^{2}Q_{2}+\frac{1}{2}%
Q_{1}^{2}a^{\prime }-\frac{1}{2}Q_{2}^{2}a^{\prime }+  \nonumber \\
&&+B_{2}Q_{1}^{\prime }+\frac{3}{4}aQ_{1}Q_{1}^{\prime }-\frac{1}{2}%
aQ_{2}Q_{2}^{\prime }+a^{\prime }Q_{2}^{\prime }+\frac{1}{2}aQ_{2}^{\prime
\prime }.
\end{eqnarray}



\end{document}